\documentclass[11pt]{amsart}

\usepackage{amsfonts}
\usepackage{mathrsfs}
\usepackage{bbm}
\usepackage{amsmath}
\usepackage{amssymb}
\usepackage{amscd}
\usepackage{amsthm}

\makeatletter
\@namedef{subjclassname@2020}{\textup{2020} Mathematics Subject Classification}
\makeatother

\numberwithin{equation}{section}

\newcommand{\di}{\displaystyle}

\newtheorem{thm}{Theorem}[section]
\newtheorem{lemma}[thm]{Lemma}
\newtheorem{prop}[thm]{Proposition}
\newtheorem{cor}[thm]{Corollary}

\newcommand{\bc}{\mathbb{C}}
\newcommand{\bn}{\mathbb{N}}
\newcommand{\bp}{\mathbb{P}}

\newcommand{\br}{\mathbb{R}}
\newcommand{\bz}{\mathbb{Z}}

\newcommand{\sd}{\mathrm{od}}
\newcommand{\scp}{\mathscr{P}}

\newcommand{\sj}{\mathscr{J}}

\newcommand{\fp}{\mathfrak{P}}

\newcommand{\jnmpo}{J^{n-1}\left(t_0, \mathbb{P}^1\right)}

\newcommand{\jnmc}{J^{n-1}\left(t_0, \mathbb{C}\right)}

\newcommand{\sjt}{\mathscr{J}_{t_0}}

\newcommand{\sjm}{\mathscr{J}^{n-1}_{t_0}}
\newcommand{\sjmx}{\mathscr{J}^{n-1}_{t_0} x}
\newcommand{\sjmy}{\mathscr{J}^{n-1}_{t_0} y}

\newcommand{\sjimx}{\mathscr{J}^{i_m}_{t_0} x}

\newcommand{\cisoc}{C^{\infty}\left(S^1, \mathbb{C} \right)}

\newcommand{\cisoca}{C^{\infty}\left(S^1, \mathbb{C}^{\ast} \right)}

\newcommand{\pn}{\mathscr{P}_{n}}
\newcommand{\pnd}{\mathscr{P}_{n, d}}

\newcommand{\pny}{\mathscr{P}_{n, 1}}

\newcommand{\pnn}{\mathscr{P}_{n, n}}

\newcommand{\pndl}{\mathscr{P}_{n, d, l}}

\newcommand{\aff}{\mathrm{Aff}(\mathbb{C})}

\newcommand{\zero}{\mathbf{0}}

\begin{document}

\title{Indecomposable representations of $C^{\infty}(S^1,\aff)$ on polynomial spaces over jet spaces}

\author{Xiaokun Wang}

\address{School of Information Science and Engineering, Shandong Univ., Qingdao, 266237,  P. R. China}

\email{wangxiaokun@sdu.edu.cn}

\author{Ning Zhang}

\address{School of Mathematical Sciences, Ocean Univ. of China, Qingdao, 266100, P. R. China}

\email{nzhang@ouc.edu.cn}



\thanks{This research was partially supported by the Natural Science Foundation of Shandong Province of China
grant ZR2023MA035.}

\begin{abstract}
  We study a class of finite dimensional complex vector spaces \(\pn\), \(n \in \bn\), consisting of polynomials on \(\jnmc\), the space of \((n-1)\) jets at a fixed point \(t_0 \in S^1\) of \(C^\infty\) functions \(S^1 \to \bc\). The loop group \(C^\infty(S^1,\aff)\) of the affine group \(\aff\) acts naturally on \(\pn\) by an indecomposable linear representation. The space \(\pn\) is canonically isomorphic to the space of holomorphic sections of a certain line bundle over the loop space of the Riemann sphere. It has an intricate structure, leading to infinite families of combinatorial identities. In this paper, we establish necessary conditions for a polynomial on \(\jnmc\) to belong to \(\pn\). Furthermore, we give partial answers to several conjectures concerning \(\pn\).
\end{abstract}

\maketitle

\pagestyle{myheadings} \markboth{\centerline{X. Wang and N. Zhang}}{\centerline{Indecomposable representations of $C^{\infty}(S^1,\aff)$}}



\section{Introduction \label{intro}}

Fix a point \(t_0 \in S^1\). For each \(n \in \bn\), the space \(\jnmc\) of \((n-1)\) jets at \(t_0\) of \(C^\infty\) maps \(S^1 \to \bc\) is isomorphic to \(\bc^n\). For \(x \in \cisoc\), we write
\[
\sjmx=(x(t_0), x^{(1)}(t_0), \cdots, x^{(n-1)}(t_0)) \in \bc^{n}.
\]
Let \(\rho_1, \rho_2\) be polynomials in $\sjmx$. We say that \(\rho_1\) and \(\rho_2\) are compatible with each other, and write \(\rho_1 \sim_n \rho_2\), if
\begin{equation} \label{compatible}
\rho_1\left( \sjmx \right)=x(t_0)^n \rho_2\left( \sjm x^{-1} \right), \hspace{2mm} x(t_0)\neq 0.
\end{equation}
Condition \eqref{compatible} implies that \(\rho_2 \sim_n \rho_1\).
Let \(\pn\) be the vector space of polynomials \(\rho_1\) on \(\jnmc\) for which there exists a polynomial \(\rho_2\) on \(\jnmc\) with \(\rho_1 \sim_n \rho_2\). Define \(\scp_0=\bc\) as the space of constants. We may regard elements of \(\scp_{n-1}\) as polynomials on \(\jnmc\) independent of \(\di x^{n-1}(t_0)\). Then \(\scp_{n-1} \subset \pn\), $n \in \bn$.

A novel phenomenon in infinite dimensional complex analysis is that the Dolbeault cohomology classes on the loop space $C^{\infty}(S^1,\bp^1)$ of the Riemann sphere depend on higher order jets of the loops \cite{z03,lz}. The space $\pn$ was first studied in \cite[Section~4]{z03}. It is canonically isomorphic to the space of holomorphic sections (equivalently, the zeroth order Dolbeault cohomology group) of a particular line bundle over $C^{\infty}(S^1,\bp^1)$.

The space $\pn$ was further explored in \cite{wz}, where it is regarded as the space of holomorphic sections of a certain line bundle over the space $\jnmpo$ of $(n-1)$-jets at $t_0\in S^1$ of $C^{\infty}$ maps $S^1\to\bp^1$, which is a finite dimensional rationally connected noncompact complex algebraic manifold. One of the discoveries in \cite{wz} is that the spaces $\pn$, $n\in\bn$, give rise to infinite families of combinatorial identities. In addition, three conjectures concerning the structure of $\pn$ were proposed.

Loop groups play an important role in quantum field theory and integrable systems, and they have deep connections with many other branches of mathematics (e.g., \cite{cl,f, g,hsw,ps}).
There is a natural linear representation of the loop group $C^{\infty}(S^1, \aff)$ on $\pn$, where $\aff$ is the affine group of the complex plane $\bc$. In this paper, we study the algebraic aspects of $\pn$. Our approach relies on the action of $C^{\infty}(S^1, \aff)$ on $\pn$, and we avoid complex analysis in our proofs.
Although the results obtained here are helpful for understanding $C^{\infty}(S^1,\bp^1)$ as a complex manifold.

Let $\scp_{n,d}$ be the subspace of $\pn$ consisting of homogeneous polynomials of degree $d$, where $d\ge 0$. Clearly, $\di\scp_{n,0}$ is the space of constants. The following results were proved in \cite[Sections 3, 4]{z03} (see also \cite[Section 2]{wz} for a detailed summary): $$\pn=\bigoplus_{d=0}^{n}\pnd;$$
for every $\rho\in\pn$, its degree as a polynomial in $\di x^{(n-1)}(t_0)$ is at most $1$; and the monomials $\di x(t_0)$, $\di x^{(1)}(t_0)$, $\dots$, $\di x^{(n-1)}(t_0)$ form a basis for $\pny$.

Given a monomial $\di\rho_{\mu}=c_{\mu}x^{(r_1)}(t_0)\cdots x^{(r_m)}(t_0)$ with $c_{\mu} \in \bc^{\ast}=\bc\setminus\{0\}$, we call $r_1+\cdots+r_m$ the order of derivatives of $\rho_{\mu}$ and denote it by $\sd(\rho_{\mu})$.
The order of derivatives of a constant function is defined to be $0$. For a nonzero polynomial $\rho=\sum_{\mu}\rho_{\mu}$ in finitely many derivatives $\di x^{(\nu)}(t_0)$, we define
\[
    \sd(\rho)=\max_{\mu}\left\{\sd(\rho_{\mu})\right\}.
\]
Let $\pndl\subset\pnd$ be the subspace spanned by polynomials $\rho\in\pnd$ such that the order of derivatives of every monomial in $\rho$ is $l$.

It was proved in \cite{wz} that
\begin{equation}\label{sym}
    \dim\pndl=\dim\scp_{n,n-d,l},\hspace{1.5mm} n\in\bn,\ d,l\in\bn\cup\{0\},\ d\le n
\end{equation}
(for more details see Section~\ref{bg}). Hence $\pnn$ is spanned by $x(t_0)^n$ and $\dim\scp_{n,n-1}=n$.
By \cite[Proposition 3.1]{wz}, a polynomial in $\sjmx$ belongs to $\pndl$ with $l+d>n$ if and only if certain combinatorial identities hold.
Using this characterization, a basis for $\scp_{n,2}$ was computed in \cite[Theorem 1.1]{wz}.
Moreover, it was proved in \cite[Proposition 3.2]{wz} that for every $n\in\bn$, there exists $N=N(n)\in\bn$ such that
\begin{equation} \label{pndl}
   \pnd=\oplus_{l=0}^{N} \scp_{n, d, l}.
\end{equation}
The following conjectures were proposed in \cite[Section 1]{wz}:
\begin{itemize}
  \item[(i)] $\di \dim \pnd=\binom{n}{d}$ for $0 \le d \le n$.

  \item[(ii)] $\di \pnd=\oplus_{l=0}^{dn-d^2} \scp_{n, d, l}$, and $\dim \scp_{n, d, l} \neq 0$ for $0 \le l \le dn-d^2$.

  \item[(iii)] $\di \dim \scp_{n, d, l}=\dim \scp_{n, d, dn-d^2-l}$ for $0 \le l \le dn-d^2$.
\end{itemize}
It was proved in \cite{z03,wz} that (i)--(iii) hold for all $n\ge1$ and $d\in\{0,1,2,n-2,n-1,n\}$; consequently, (i)--(iii) are valid for $1\le n\le5$.

Let $\rho$ be a polynomial in $\sjmx$. For $m=1,\dots,n$, let $M_{m,\rho}$ be the set of products
\begin{equation}\label{product}
    x^{(r_1)}(t_0)\cdots x^{(r_m)}(t_0)
\end{equation}
such that at least one monomial in $\rho$ is divisible by (\ref{product}).
Define
\begin{equation} \label{odm}
\sd_{m}(\rho)=\max_{x^{(r_1)}(t_0) \cdots x^{(r_m)}(t_0) \in M_{m, \rho}}(r_1+\cdots+r_m).
\end{equation}
If $M_{m, \rho}=\emptyset$, we set $\sd_{m}(\rho)=0$. Clearly, $\sd_m(\rho) \le \sd(\rho)$ for $m=1,\dots,n$. Moreover, if $\rho \in \pnd$, then $\sd_d(\rho)=\sd(\rho)$.

The main results of this paper are as follows.
\begin{thm} \label{arb}
Let $\rho \in \pn$. Then
\begin{equation} \label{upperb}
\sd_m(\rho) \le m(n-m), \hspace{1.5mm} m=1,\dots,n.
\end{equation}
In particular, $\di \pnd=\oplus_{l=0}^{dn-d^2} \scp_{n, d, l}$.
\end{thm}

By Theorem \ref{arb}, the first part of conjecture (ii) is answered affirmatively.
Note that the estimate (\ref{upperb}) is sharp for $m=0,1,2,n-2,n-1,n$.
It turns out that if $n \ge 3$, $2 \le d \le (n+1)/2$, $\dim \scp_{n, d, dn-d^2} \neq 0$ and $\rho \in \scp_{n, d, dn-d^2} \setminus \{\zero\}$, then $\sd_m(\rho)=m(n-m)$ for every $m=1,\dots,d$ (Corollary \ref{2}).

We already know that if $2 \le d \le n-2$, then $\dim \scp_{n, d, 0} =1$, $\dim \scp_{n, d, 1} =1$ and $\dim \scp_{n, d, 2} =2$ (Proposition \ref{monoj}).
As another application of Theorem \ref{arb}, we obtain the following partial answer to conjecture (iii).

\begin{thm} \label{bounds}
The following statements hold for $n \ge 6$ and $3 \le d \le n-3$.
\begin{itemize}
  \item[(a)] $\dim \scp_{n, d, dn-d^2-l} \le \dim \scp_{n, d, l}$ for $l=0,1,2$.

  \item[(b)] If $n-2d+1 \ge 0$ and $\dim \scp_{n, d, dn-d^2} \neq 0$, then for $0 \le l \le d(n-2d+1)$,
  $\dim \scp_{n, d, dn-d^2-l}$ is at least the number of elements in the set
  $S_{n, d, l}$ of multi-indices
  $(k_1,\dots,k_d) \in \bz^d$ such that $0 \le k_1 \le \cdots \le k_d \le n-2d+1$ and $\sum_{i=1}^d k_i=l$.
  In particular, if $l \le n-2d+1$, then the number of elements in $S_{n, d, l}$ is precisely $\dim \pndl$.
\end{itemize}
\end{thm}

This paper is organized as follows.  In Section~\ref{bg}, we recall some relevant results from \cite{z03,wz}. In Section~\ref{action}, we describe a natural representation of $C^{\infty}(S^1,\aff)$ on $\pn$ and show that this representation is indecomposable; we also study relations between different elements of $\pn$. In Section~\ref{proof}, we prove Theorem~\ref{arb}. In the final Section~\ref{1.2}, we prove Theorem~\ref{bounds}. The space $\scp_{n,d,dn-d^2}$ plays a key role in the proof of Theorem~\ref{bounds}(b). The proofs of both Theorem~\ref{arb} and Theorem~\ref{bounds} use the results of Section~3.


\section{Background \label{bg}}

The space $\pn$ admits three equivalent definitions. It was first introduced on P. 12 of \cite{z03}
as the space $\mathrm{R}(\mathfrak{H}_{\varphi})$ (with $r=1$, $t_1=t_0$ and
$n_1=n$ in \cite[(4.1)]{z03}), consisting of
local representations of holomorphic sections of a particular line bundle
on the loop space of the Riemann sphere.
It was subsequently defined as the space of local representations of holomorphic sections of a particular
line bundle on the space $\jnmpo$ of $n-1$ jets of $C^{\infty}$
maps $S^1 \to \bp^1$ at a given point $t_0 \in S^1$ (P. 591 of \cite{wz}).
These two definitions are equivalent to the one given at the beginning of this paper (\cite[Section 2]{wz}).

A complex manifold is said to be rationally connected if any two points can be joined by a chain of rational curves, i.e., by a chain of holomorphic images of $\mathbb{P}^1$.
Rationally connected algebraic manifolds are a highly active research area in algebraic geometry.
Recall that the loop space of the Riemann sphere $C^{\infty}(S^1, \bp^1)$ is rationally connected
(\cite[Proposition 3.3]{z17}, where $\mathscr{C}$ was defined on P. 711 of \cite{z17}).
There is a holomorphic surjective map $C^{\infty}(S^1, \bp^1) \to \jnmpo$ sending each $x \in C^{\infty}(S^1, \bp^1)$
to its $(n-1)$ jet at $t_0 \in S^1$ (P. 589 of \cite{wz}). Hence $\jnmpo$
is rationally connected.

Let $\rho\in\pn$. Suppose that $x^{(k_0)}(t_0)$ is the highest derivative on which $\rho$ depends, and let $m$ be the degree of $\rho$ as a polynomial in $x^{(k_0)}(t_0)$. It follows from \cite[Proposition 4.3]{z03} that
\begin{equation}\label{highest}
    m(k_0+1)\le n.
\end{equation}
Hence $k_0\le n-1$. Moreover, if $k_0>n/2-1$, then $m=1$.

For $k \in \bn$, let $\fp_{k}$ denote the space of $k$-tuples of nonnegative integers
$\di p=(p_1, \cdots, p_k)$ satisfying $$\di p_1+2p_2+\cdots +k p_k=k,$$ and set $|p|=p_1+\cdots+p_k$.
Given $\di x \in \cisoc$ with $x(t_0) \not=0$, the Fa\`{a} di Bruno's formula \cite{hu} yields
\begin{eqnarray} \label{chain}
& \di \left(\frac{1}{x} \right)^{(k)}(t_0)  =  \sum_{p \in \fp_k}  c_p \frac{x^{(1)}(t_0)^{p_{1}}
\cdots x^{(k)}(t_0)^{p_{k}}}{x(t_0)^{|p|+1}}, \hspace{1.5mm} \text{where} & \\
& \di
c_p=k ! \frac{(-1)^{|p|}}{(1!)^{p_1}  \cdots (k!)^{p_k}} \binom{|p|}{p_1, \cdots, p_k} \hspace{2mm}
\text{and} \hspace{2mm} \binom{|p|}{p_1, \cdots, p_k}=\frac{|p|!}{p_1! \cdots p_k!}. & \notag
\end{eqnarray}
On the right hand side of (\ref{chain}), we have $1 \le |p| \le k$.
The only term with $|p|=1$ respectively $|p|=k$ is
\begin{equation} \label{p=1}
\di \frac{-x^{(k)}(t_0)}{x(t_0)^2}
\hspace{3mm} \text{respectively} \hspace{3mm}  \frac{(-1)^k k! x^{(1)}(t_0)^k}{x(t_0)^{k+1}},
\end{equation}
and the order of derivatives of the numerator of every term is $k$.

Let $\rho_{\mu} \not=0$ be a monomial in $\sjmx$. It follows from (\ref{chain}) that
the function $\di \rho_{\mu} \left( \sj_{t_0}^{n-1} x^{-1} \right)$, $x(t_0) \not=0$, can be written
as a finite linear combination of pairwise distinct rational functions in
$x(t_0)$, $x^{(1)}(t_0)$, $\cdots$, $x^{(n-1)}(t_0)$. In each rational function,
the numerator is
a monic monomial in $x^{(1)}(t_0)$, $\cdots$, $x^{(n-1)}(t_0)$ whose order of derivatives is  $\sd(\rho_{\mu})$, and
the denominator is a power of $x(t_0)$.
We call this linear combination the inversion expansion of $\di \rho_{\mu}$.
For a polynomial $\rho =\sum_{\mu} \rho_{\mu}$, the inversion expansion is defined
as the sum of the inversion expansions of its monomials $\rho_{\mu}$.
The right hand side of (\ref{compatible}) may be viewed as the product of $x(t_0)^n$
and the inversion expansion of $\rho_2$. Thus
\begin{equation} \label{sd}
\sd(\rho_1)=\sd(\rho_2) \hspace{2mm} \text{if} \hspace{2mm} \rho_1 \sim_n \rho_2.
\end{equation}
We shall repeatedly use the following result
\cite[Proposition 3.1]{wz}:

\begin{prop} \label{monoj}
    Let $\rho$ be a polynomial in $x(t_0)$, $x^{(1)}(t_0)$, $\cdots$, $x^{(m)}(t_0)$.
    Then $\rho \in \scp_{n} \setminus \scp_{n-1}$ if and only if,
    in the inversion expansion of $\rho$, the highest power of $x(t_0)$  appearing in the denominators is $x(t_0)^{n}$.
    In particular, the monomial
    \begin{equation*}
  x^{(k_1)}(t_0) \cdots x^{(k_d)}(t_0) \in \scp_{k_1+\cdots+k_d+d} \setminus \scp_{k_1+\cdots+k_d+d-1},
    \end{equation*}
    where $d=1, 2, \cdots$ and $k_1, \cdots, k_d$ are nonnegative integers.
\end{prop}

By results on P. 594 of \cite{wz} and (\ref{sd}), we have  the following
\begin{prop} \label{pnd}
  Let $\rho_1 \in \pnd$ and $\rho_2 \in \pndl$. If $\rho_1 \sim_n \rho_3$ and $\rho_2 \sim_n \rho_4$,
  then $\rho_3 \in \scp_{n, n-d}$ and $\rho_4 \in \scp_{n, n-d, l}$. Thus $\di \dim \pnd=\dim \scp_{n, n-d}$
  and $\di \dim \pndl=\dim \scp_{n, n-d, l}$.
\end{prop}

\section{The action of $C^{\infty}(S^1, \aff)$ on $\pn$ \label{action}}

The main result of this section is Theorem \ref{monomial}.

\begin{lemma} \label{invariantr}
Let $\di\rho(\sjmx)\in\pn$. Fix $\di y_0,z_0\in\cisoc$ with $y_0(t_0)\neq0$. Then $\rho(\sjm(y_0x+z_0))$, regarded as a polynomial in $\sjmx$, belongs to $\pn$ and has the same order of derivatives as $\rho(\sjmx)$.
\end{lemma}
\begin{proof}
  (i) Replacing $x$ by $y_0x$ in (\ref{compatible}), with $(y_0x)(t_0) \not=0$, we obtain
  \begin{eqnarray*}
  & \rho_1( \sjm(y_0x))=x(t_0)^n \rho_3(\sjm x^{-1}), \hspace{1mm} \text{where} & \\
  & \rho_3(\sjm x^{-1})= y_0(t_0)^n \rho_2( \sjm (x^{-1} y_0^{-1})). &
  \end{eqnarray*}
  It is clear that
  $y_0(t_0)^n \rho_2( \sjm (x^{-1} y_0^{-1}))$ is a polynomial in $\sjm x^{-1}$,
  whose coefficients depend on the constants $y_0(t_0)$, $y^{(1)}_0(t_0)$, $\cdots$, $y_0^{(n-1)}(t_0)$. By definition,
  $\rho_1(\sjm(y_0x)) \in \pn$.

  (ii) Let $\zeta=\zeta(z_0, x)=1+z_0x^{-1}$, where $x(t_0) \not=0$. Then $\zeta$ is a well defined
  $C^{\infty}$ function in a neighborhood of $t_0 \in S^1$.
  Replacing $x$ by $x+z_0$ in (\ref{compatible}),
  with $x(t_0) \not=0$ and $x(t_0)+z_0(t_0) \not=0$ (hence
 $\zeta(t_0) \not=0$),
we obtain
  \begin{eqnarray} \label{ga}
 & \rho_1(\sjm(x+z_0))=x(t_0)^n  \rho_4\left( \sjm x^{-1} \right), \hspace{1mm} \text{where} & \\
 & \rho_4\left( \sjm x^{-1} \right)= \zeta(t_0)^n \rho_2(\sjm(x^{-1} \zeta^{-1})). & \notag
  \end{eqnarray}
 It suffices to show that $ \rho_4\left( \sjm x^{-1} \right)$ is a polynomial in $\sjm x^{-1}$.
  By the Leibniz rule
  and (\ref{chain}), we can express
  $\rho_2(\sjm(x^{-1} \zeta^{-1}))$ as
  a finite linear combination of pairwise distinct rational functions,
  where the denominator of each rational function is a power of $\zeta(t_0)$ and the numerator
  is a monic monomial in $(x^{-1})^{(l)}(t_0)$, $l=0, 1, \cdots, n-1$, and
  $\zeta^{(i)}(t_0)$, $i=1, \cdots$, $n-1$.
  In turn, each $\zeta^{(i)}(t_0)$, $i=0, 1, \cdots$, $n-1$, is a polynomial in $\sjm x^{-1}$.
  Suppose that the highest power of $\zeta(t_0)$ in the denominators is $\zeta(t_0)^{N}$.

  Assume for contradiction that $N>n$.
  Let $z_{\lambda}=z_0+\lambda x \in \cisoc$, where $\lambda \in \bc$ and $x(t_0) \not=0$.
  Write $\zeta_{\lambda}$ for $1+z_{\lambda} x^{-1}=\zeta+\lambda$ (which is
  well defined in a neighborhood of $t_0 \in S^1$), so that $\zeta_0=\zeta$.
  Note that $\zeta_{\lambda}^{(i)}(t_0)=\zeta^{(i)}(t_0)$ for $i=1, \cdots, n-1$. Now
  replace $z_0$ by $z_{\lambda}$ in (\ref{ga})
  (thus $\zeta$ is replaced by $\zeta_{\lambda}$), multiply both sides
  of the equation by $\zeta_{\lambda}(t_0)^{N-n}$, and let $\lambda \to -\zeta(t_0)$.
  Then $\zeta_{\lambda}(t_0)\to0$, so for any given $x\in\cisoc$ the left hand side of the equation converges to $0$,
  whereas the right hand side does not.
  This is a contradiction.
  Hence $N \le n$, and $ \rho_4\left( \sjm x^{-1} \right)$ is a polynomial in $\sjm x^{-1}$.

  (iii) Composing the transformations in (i) and (ii), we conclude that $\rho(\sjm(y_0x+z_0)) \in \pn$. A direct verification shows that its order of derivatives is the same as that of $\di \rho(\sjm x)$.
\end{proof}

Let $\aff$ be the group of affine transformations on $\bc$, which is precisely
the group of holomorphic automorphisms of $\bc$.
The group $\aff$ is isomorphic to the semidirect product $\bc^{\ast} \ltimes \bc$,
where $\bc^{\ast}$ acts on $\bc$ by multiplication.
Consequently, the loop group
$$C^{\infty}(S^1, \aff) \simeq \cisoca \ltimes \cisoc,$$
where $\cisoca$ acts on $\cisoc$ by multiplication, and
the group operation is given by
$$(y_0, z_0)(y_1, z_1)=(y_0 y_1, y_1z_0+z_1).$$
By Lemma \ref{invariantr}, there exists a left linear representation
$$\pi: \cisoca \ltimes \cisoc \times \pn \to \pn, $$ defined by
\begin{equation*}
 \pi((y, z), \rho(\sjm x))=\rho( \sjm(yx+z)).
\end{equation*}
We denote this representation by $(\pi, \pn)$.

\begin{prop}  \label{prepare}
  Let $(\pi, V)$ be a subrepresentation of $(\pi, \pn)$.
\begin{enumerate}
  \item[(a)]
  Let
  $\rho=\sum_{j=1}^m \rho_{d_j} \in V \setminus \{\mathbf{0}\}$,
  where $n \ge d_1>\cdots >d_m \ge 0$ and
  $\rho_{d_j} \in \scp_{n, d_j} \setminus \{\mathbf{0}\}$.
  Then $\rho_{d_j} \in V$ for every $j=1,\cdots,m$. Consequently,
  \begin{equation*}
    V=\bigoplus_{d=0}^{n} (V \cap \pnd).
  \end{equation*}

  \item[(b)]
  Let $\tilde{\rho} \in V$ and $x, y \in \cisoc$. We may write
  \begin{eqnarray*}
    \tilde{\rho}(\sjm (x+y))
    &=& \tilde{\rho}(\sjmy) \\
    &+& \sum_{d=1}^n \sum
    \tilde{\rho}_{(k_1, \cdots, k_d)}(\sjm y)
    x^{(k_1)}(t_0) \cdots x^{(k_d)}(t_0),
  \end{eqnarray*}
  where the functions $\tilde{\rho}_{(k_1, \cdots, k_d)}(\sjm y)$ are polynomials in
  $\sjm y$, and the second summation is taken over all multi-indices
  $(k_1$,$\cdots$, $k_d)$ $\in \bz^d$
  such that
  $k_1 \ge \cdots \ge k_d \ge 0$ and at least one of the monomials in
  $\tilde{\rho}(\sjmx)$ is divisible by
  $x^{(k_1)}(t_0) \cdots x^{(k_d)}(t_0)$.
  Then every polynomial $\tilde{\rho}_{(k_1, \cdots, k_d)}(\sjm y)$ belongs to $V$.

  \item[(c)]
  Let
  $\rho^{\ast}=\sum_{d=0}^m \rho^{\ast}_{d} \in V$,
  where $\rho^{\ast}_d \in \pnd \cap V$, and let $x \in \cisoc$ and
  $y \in \cisoca$. We may write
  \begin{eqnarray*}
    \rho^{\ast}(\sjm (xy))
    &=& \rho^{\ast}_0(\sjmy) \notag \\
    &+& \sum_{d=1}^n \sum
    \rho^{\ast}_{(k_1,\cdots,k_d)}(\sjm y)
    x^{(k_1)}(t_0) \cdots x^{(k_d)}(t_0), \label{xy}
  \end{eqnarray*}
  where the functions $\rho^{\ast}_{(k_1,\cdots,k_d)}(\sjm y)$ are polynomials
  in $\sjm y$, and the second summation is taken over all multi-indices
  $(k_1$, $\cdots$, $k_d)$ $\in \bz^d$ such that
  $k_1 \ge \cdots \ge k_d \ge 0$ and at least one of the monomials in
  $\rho^{\ast}_d(\sjmx)$ is divisible by a product
  $x^{(r_1)}(t_0) \cdots x^{(r_d)}(t_0)$ with $r_1 \ge \cdots \ge r_d \ge 0$ and
  $k_j \le r_j$ for $j=1,\cdots,d$.
  Then every polynomial $\rho^{\ast}_{(k_1,\cdots,k_d)}(\sjm y)$ belongs to
  $\pnd \cap V$.
\end{enumerate}
\end{prop}
\begin{proof}
(a) Note that
  $\rho(\sjm (\lambda x))
    =\sum_{j=1}^m \lambda^{d_j} \rho_{d_j}(x) \in V$,
  where $\lambda \in \bc^{\ast}$ and $x \in \cisoc$. The polynomials
  $\lambda^{d_j}$, $j=1,\cdots,m$, on $\bc$ are linearly independent.
  Hence there exist pairwise distinct points
  $\lambda_1,\cdots,\lambda_m \in \bc^{\ast}$ such that the
  $m \times m$ matrix
  $(\lambda_i^{d_j})_{i,j=1,\cdots,m}$ is invertible.
  Let
  $g_i=\sum_{j=1}^m \lambda^{d_j}_i \rho_{d_j} \in V$,
    $i=1,\cdots, m$.
  Then each $\rho_{d_j}$ is a finite linear combination of
  $g_1,\cdots,g_m$. Therefore $\rho_{d_j} \in V$ for every
  $j=1,\cdots,m$.

  (b), (c) The monomials
  $x^{(k_1)}(t_0) \cdots x^{(k_d)}(t_0)$ in (b) or (c) and the constant $1$ are linearly
  independent. The conclusions of (b) and (c) follow from arguments
  similar to those in (a).
\end{proof}

For a general multi-index $(k_1,\cdots,k_d) \in \bz^d$ with
$k_1 \ge \cdots \ge k_d \ge 0$, the product
$x^{(r_1)}(t_0) \cdots x^{(r_d)}(t_0)$ appearing in
Proposition \ref{prepare}(c) need not exist. In that case, we simply
set $\rho^{\ast}_{(k_1,\cdots,k_d)}=\zero$. Consequently, there is always
a well defined linear operator
\begin{equation} \label{mul}
  \mathfrak{M}_{(k_1,\cdots,k_d)}: V \ni \rho^{\ast}
  \mapsto \rho^{\ast}_{(k_1,\cdots,k_d)} \in \pnd \cap V.
\end{equation}
Note that if $\rho^{\ast} \in \pndl \cap V$, then $\mathfrak{M}_{(k_1,\cdots,k_d)}(\rho^{\ast}) \in \scp_{n, d, l-k_1-\cdots-k_d}$.

Let $\rho(\sjm x) \in \pn$. Suppose that at least one of the monomials in $\rho$ is divisible by the product
$$
\di x^{(i_1)}(t_0) \cdots x^{(i_m)}(t_0), \hspace{1mm}  m \ge 1, \hspace{1mm}  i_1  \ge \cdots \ge i_m \ge 0.
$$
Then we may write
\begin{equation} \label{decomp}
\rho(\sjm x)= x^{(i_1)}(t_0) \cdots x^{(i_m)}(t_0) f_1(\sjm x)+f_2(\sjm x),
\end{equation}
where $f_1(\sjm x)$ and $f_2(\sjm x)$ are polynomials, $f_1 \neq \zero$, and, if $f_2 \neq \zero$, no monomial in $f_2(\sjm x)$ is divisible by $\di x^{(i_1)}(t_0) \cdots x^{(i_m)}(t_0)$.
Recall the notation $\sd_m(\rho)$ as defined in (\ref{odm}).

\begin{prop} \label{1}
   For any subrepresentation $(\pi,V)$ of $(\pi,\pn)$, the constant function $1$ belongs to $V$.
\end{prop}
\begin{proof}
  Assume $V \neq \scp_0 = \bc$. Choose a nonconstant polynomial $\rho_0 \in V$. Let $x^{(N_0)}(t_0)$ be the highest derivative on which $\rho_0$ depends, and let $m_0 \ge 1$ be the degree of $\rho_0$ as a polynomial in $x^{(N_0)}(t_0)$. Express $\rho_0$ in the form (\ref{decomp}) by setting $m=m_0$ and $i_1=\cdots=i_{m_0}=N_0$. Then the polynomial $f_1 \not=\zero$ in (\ref{decomp}) is independent of $x^{(k)}(t_0)$ for $N_0 \le k \le n-1$.

Consider $\di \rho(\sjm(x+y_0))$ as a polynomial in $\di \sjmx$. One of its monomials is of the form
$f_1(\sjt^{N_0-1}y_0) x^{(N_0)}(t_0)^{m_0}$.
By Proposition \ref{prepare}(b), we have $f_1(\sjt^{N_0-1} y_0) \in V \setminus \{\zero\}$.
If $N_0=0$, then $f_1$ is constant, so $1 \in V$.
Otherwise, set $\rho_1=f_1$, and let $x^{(N_1)}(t_0)$ be the highest derivative on which $\rho_1$ depends. Then $N_1<N_0$. Repeating the above process, we obtain $1 \in V$ after finitely many steps.
\end{proof}

We say that a subrepresentation $(\pi, V)$ of $(\pi, \pn)$ is trivial if $V=\{\mathbf{0}\}$ or $V=\pn$.
Recall that $(\pi, \scp_{n-1})$ is a nontrivial subrepresentation of $(\pi, \pn)$ for $n \in \bn$.
As an immediate consequence of Proposition \ref{1}, we obtain the following

\begin{cor} \label{inde}
  The representation $(\pi, \pn)$ is indecomposable; that is, it is not the direct sum of two nontrivial subrepresentations.
\end{cor}

\begin{thm}
  \label{monomial}
  Let $\rho \in \pn$ be a polynomial of degree $d=1, \cdots, n$.
  Then for any $m=1, \cdots, d$, there exists $\tilde{\rho}_m \in \scp_{n, m, \sd_m(\rho)}$ of the form
    \begin{equation} \label{m}
    \sum \tilde{c}_{(j_1, \cdots, j_m)}x^{(j_1)}(t_0) \cdots x^{(j_m)}(t_0),
    \end{equation}
    where every
    coefficient $\tilde{c}_{(j_1, \cdots, j_m)} \in \bc^{\ast}$,
    and
    the summation ranges over all multi-indices $(j_1, \cdots, j_m)$ such that
    $j_1  \ge \cdots \ge j_m \ge 0$, $j_1+\cdots+j_m=\sd_m(\rho)$
     and at least one of the monomials in $\rho$ is divisible by
    the product $x^{(j_1)}(t_0) \cdots x^{(j_m)}(t_0)$.
\end{thm}
\begin{proof}
Choose the product
\(x^{(i_1)}(t_0)\cdots x^{(i_m)}(t_0)\)
in \eqref{decomp} such that
\begin{equation}\label{maximal}
i_1  \ge \cdots \ge i_m \ge 0 \hspace{2mm} \text{and} \hspace{2mm}  i_1+\cdots+i_m=\sd_m(\rho).
\end{equation}
We claim that the polynomial \(f_1\) in \eqref{decomp} does not involve \(x^{(k)}(t_0)\) for \(i_m+1\le k\le n-1\).
Indeed, otherwise we could replace \(x^{(i_m)}(t_0)\) by a higher order derivative \(x^{(k)}(t_0)\) occurring in \(f_1\). This would yield a product of \(m\) derivatives of \(x\) at \(t_0\) whose order of
 derivatives is strictly larger than \(\sd_m(\rho)\). Since this product still divides a monomial of \(\rho\), it contradicts the maximality in \eqref{maximal}.

  Suppose that
\[
i_1\ge\cdots\ge i_{m-I}>i_{m-I+1}=\cdots=i_m,
\]
where \(1\le I\le m\) (omit \(i_{m-I}\) when \(I=m\)).
Let \(J\ge 0\) be the degree of \(f_1\) as a polynomial in \(x^{(i_m)}(t_0)\).
Then we may write
\(x^{(i_1)}(t_0)\cdots x^{(i_m)}(t_0)\, f_1(\sjimx)\)
in the form
  \begin{equation} \label{gpaction}
  x^{(i_1)}(t_0) \cdots x^{(i_{m-I})}(t_0) \sum_{s=I}^{I+J}   g_s(\sjt^{i_m-1}x) x^{(i_m)}(t_0)^s,
  \end{equation}
  where each \(g_s\) is a polynomial in $\sjt^{i_m-1} x$ (all \(g_s\) are constants if \(i_m=0\)), and \(g_{I+J}\neq \zero\).

  For any fixed $\di y_0 \in \cisoca$,
  replace $x$ by $y_0x+y_0$ in (\ref{decomp}). It follows from Lemma \ref{invariantr} that
  $\di \rho(\sjm(y_0x+y_0))$,  viewed as a polynomial in $\di \sjmx$, belongs to $\pn$.
  For each \(l=1,\dots,n-1\),
  regard $(y_0x+y_0)^{(l)}(t_0)$
  as a polynomial in $\sjt^{l} x$; the
   only monomial involving $x^{(l)}(t_0)$ is then $y_0(t_0) x^{(l)}(t_0)$.
  Recall that no monomial of $f_2(\sjm x)$ is divisible
by $\di x^{(i_1)}(t_0) \cdots x^{(i_m)}(t_0)$ whenever $f_2 \not= \zero$.
  By (\ref{maximal}), no
  monomial in $\di f_2(\sjm(y_0x+y_0))$ is divisible
  by $\di x^{(i_1)}(t_0) \cdots x^{(i_m)}(t_0)$.
  Consider the monomial of the form
  $$c_{(i_1, \cdots, i_m)}(\sjt^{i_m} y_0) x^{(i_1)}(t_0) \cdots x^{(i_m)}(t_0)$$
  in
  $\di \rho(\sjm(y_0x+y_0))$.
  It follows from (\ref{decomp}) and (\ref{gpaction}) that
  \begin{equation*}
  c_{(i_1, \cdots, i_m)}(\sjt^{i_m} y_0)=y_0(t_0)^m \sum_{s=I}^{I+J} \binom{s}{I} y_0^{(i_m)}(t_0)^{s-I} g_s(\sjt^{i_m-1}y_0),
  \end{equation*}
  where $\binom{s}{I} y_0^{(i_m)}(t_0)^{s-I} y_0(t_0)^I x^{(i_m)}(t_0)^I$ is one of the monomials in
  $$\di (y_0x+y_0)^{(i_m)}(t_0)^s.$$
  Recall that the polynomial $\di g_{I+J} \not=\zero$. Hence the polynomial $c_{(i_1, \cdots, i_m)} \not=\zero$.
  Therefore $c_{(i_1, \cdots, i_m)}(\sjt^{i_m} y_0)\not=0$ whenever $\sjt^{i_m} y_0$
  is contained in a dense open subset of $J^{i_m}(t_0, \bc)$.

  Similarly,
  for every multi-index $(j_1, \cdots, j_m)$ in (\ref{m}),
  consider the monomial
  \begin{equation*} \label{c}
  c_{(j_1, \cdots, j_m)}(\sjt^{j_m} y_0) x^{(j_1)}(t_0) \cdots x^{(j_m)}(t_0)
  \end{equation*}
  occurring in
  $\di \rho(\sjm(y_0x+y_0))$.
Choose $y_0^{\ast} \in \cisoca$ such that every polynomial $c_{(j_1, \cdots, j_m)}$
is nonzero at $\sjt^{i_m} y_0^{\ast}$.
Given a multi-index $(j'_1,\dots,j'_m)$ with $j'_1\ge\cdots\ge j'_m\ge0$ and $j'_1+\cdots+j'_m=\sd_m(\rho)$, if no monomial of $\rho$ is divisible by $x^{(j'_1)}(t_0)\cdots x^{(j'_m)}(t_0)$, then no monomial of $\di\rho(\sjm(y_0x+y_0))$ is a constant multiple of such a product.
  Let $\tilde{\rho}_m$ denote the sum
  of all monomials in the $m$-th order homogeneous component of $\di \rho(\sjm(y_0^{\ast}x+y_0^{\ast}))$
  whose order of derivatives equals $\sd_m(\rho)$ (see (\ref{pndl})).
\end{proof}

\section{Proof of Theorem \ref{arb} \label{proof}}

Recall the inversion expansion of a polynomial $\rho$ in $\sjmx$, as defined in Section~\ref{bg}.

\begin{lemma} \label{estimate}
  Let $\rho\in\pnd$. Suppose that at least one monomial in $\rho$ is divisible by the product
$x^{(i_1)}(t_0)\cdots x^{(i_m)}(t_0)$,
where $1\le m\le d$ and $0<i_m\le\cdots\le i_1$; consequently, $\rho$ can be written in the form (\ref{decomp}).
Let $f_1$ and $f_2$ be the polynomials appearing in (\ref{decomp}), with $f_1$ homogeneous of degree $d-m$.
  If $f_1$ is independent of
$x^{(k)}(t_0)$ for $i_m \le k \le n-1$, either $f_2=\zero$ or the numerator of every
rational function in the
inversion expansion of $f_2$ is not divisible
by $\di x^{(i_1)}(t_0) \cdots x^{(i_m)}(t_0)$,
  then $m+d \le n$  and $f_1 \in \scp_{n-2m, d-m}$.
\end{lemma}
\begin{proof}
  The inversion
  expansion of $\rho$
 takes the following form:
  \begin{equation} \label{ier}
  (-1)^m \frac{x^{(i_1)}(t_0) \cdots x^{(i_m)}(t_0)}{x(t_0)^{2m}} f_1\left(\sjt^{i_m-1} x^{-1} \right)+R(\sjmx),
  \end{equation}
  where $\di f_1( \sjt^{i_m-1}x^{-1})$ is the inversion expansion of $f_1$,
  and $R(\sjmx)$ is a finite linear combination of pairwise distinct rational functions in $\sjmx$,
  each with denominator a power of $x(t_0)$ and numerator a monic monomial in $x^{(1)}(t_0)$, $\cdots$, $x^{(n-1)}(t_0)$.
  The conditions in the lemma imply that no numerator of a rational function in $R(\sjmx)$ is divisible by $\di x^{(i_1)}(t_0)\cdots x^{(i_m)}(t_0)$.
  By (\ref{compatible}) and Proposition \ref{pnd}, the product of $x(t_0)^n$ with the expression in (\ref{ier})
  is an element of $\scp_{n,n-d}$.
  Hence
  $$(-1)^m x^{(i_1)}(t_0) \cdots x^{(i_m)}(t_0) (x(t_0)^{n-2m} f_1(\sjt^{i_m-1} x^{-1}))$$
  is a homogeneous polynomial of degree $n-d$.
  Thus $m+d\le n$, and
  $$x(t_0)^{n-2m}f_1(\sjt^{i_m-1}x^{-1})$$
  is a homogeneous polynomial of degree $n-d-m$ in $\sjt^{i_m-1}x$.
  By Proposition \ref{monoj}, we have $f_1\in\scp_{n-2m,d-m}$.
  \end{proof}

\begin{prop} \label{homo}
  For any $\rho\in\pnd$, we have $\di\sd(\rho)\le d(n-d)$.
\end{prop}
\begin{proof}
  Without loss of generality, we may assume that $\rho\in\pndl$, where $l=\sd(\rho)$ and $d\le n/2$, see (\ref{pndl}) and Proposition~\ref{pnd}. The conclusion of the proposition is already known to hold for $d=0,1,2$. Given $d_0$ with $3\le d_0\le n/2$, suppose that the conclusion holds for all $0\le d\le d_0-1$. We then show that it also holds for $d=d_0$.

  Let $\di x^{(N)}(t_0)$ be the highest derivative on which $\rho$ depends, and let $m$ with $1\le m\le d_0$ be the degree of $\rho$ as a polynomial in $\di x^{(N)}(t_0)$. We may assume that $N>0$; otherwise $\sd(\rho)=0$. Write
\[
    \rho=x^{(N)}(t_0)^m f_1+f_2,
\]
where $f_1 \not= \zero$ is independent of $\di x^{(N)}(t_0)$ and $f_2$ has degree less than $m$ as a polynomial in $\di x^{(N)}(t_0)$.

If $m=d_0$, then $f_1$ is constant. Since all monomials of $\rho$ have the same order of derivatives,
$\rho$ consists of a single monomial; hence $\rho=c x^{(N)}(t_0)^{d_0}$ for some $c \in \bc^{\ast}$. By Proposition~\ref{monoj}, this monomial lies in $\scp_{n,d_0}$ if and only if $\di d_0N+d_0\le n$.
    Thus $$\di \sd(\rho)=d_0 N \le n-d_0 < d_0(n-d_0).$$

    If $m<d_0$, it follows from Lemma \ref{estimate} that $\di f_1 \in \scp_{n-2m, d_0-m}$.
    By the induction hypothesis, we have
    $$\sd(f_1) \le (d_0-m)(n-2m-(d_0-m))=nd_0-nm+m^2-d_0^2.$$
    By Theorem \ref{monomial} and the induction hypothesis, we obtain
    $$\sd( x^{(N)}(t_0)^m) \le m(n-m).$$
    Therefore
    $\sd(\rho)=\sd( x^{(N)}(t_0)^m)+\sd(f_1) \le d_0(n-d_0).$
\end{proof}

\noindent {\it Proof of Theorem \ref{arb}.}
Let $d=0,1,\dots,n$ be the degree of $\rho$. If $d<m\le n$, then $\sd_m(\rho)=0$, and (\ref{upperb}) follows. If $1\le m\le d$, consider the polynomial $\tilde{\rho}_m\in\scp_{n,m,\sd_m(\rho)}\setminus\{\zero\}$ in Theorem~\ref{monomial}. Applying Proposition~\ref{homo} to $\tilde{\rho}_m$, we obtain (\ref{upperb}). It follows from (\ref{pndl}) that
$\pnd=\bigoplus_{l=0}^{dn-d^2}\scp_{n,d,l}$.
\qed


\section{Upper/lower bounds for $\di \dim \pndl$ \label{1.2}}

\begin{lemma} \label{2factors}
Let $\rho\in\pn \setminus \{\zero\}$ with $n\ge2$. If $x^{(k_0)}(t_0)$ is the highest derivative on which $\rho$ depends and $2k_0>n-1$, then no monomial in $\rho$ is divisible by $x^{(k_0)}(t_0)x^{(k_0-1)}(t_0)$.
\end{lemma}
\begin{proof}
Let $m$ be the degree of $\rho$ as a polynomial in $x^{(k_0)}(t_0)$. It follows from (\ref{highest}) and the condition $2k_0>n-1$ that $m=1$. Assume for contradiction that at least one monomial in $\rho$ is divisible by $x^{(k_0)}(t_0)x^{(k_0-1)}(t_0)$.
Then $\sd_2(\rho)=2k_0-1$. By Theorem~\ref{monomial}, there exists a homogeneous polynomial $\tilde{\rho}_2 \in \scp_{n,2,2k_0-1}$ such that one monomial in $\tilde{\rho}_2$ is a nonzero constant multiple of $x^{(k_0)}(t_0)x^{(k_0-1)}(t_0)$. For any other monomial $\tilde{\rho}_\nu$ in $\tilde{\rho}_2$, there is a monomial in $\rho$ divisible by $\tilde{\rho}_\nu$.
Since $\tilde{\rho}_\nu$ has order of derivatives $2k_0-1$, it must itself be a constant multiple of $x^{(k_0)}(t_0)x^{(k_0-1)}(t_0)$.
Hence $\tilde{\rho}_2=c\,x^{(k_0)}(t_0)x^{(k_0-1)}(t_0)$ for some $c\in\bc^{\ast}$, which contradicts Proposition~\ref{monoj}.
\end{proof}

\begin{lemma} \label{n-1}
  If $\dim\scp_{n,d,dn-d^2-l}\neq0$ and $\rho\in\scp_{n,d,dn-d^2-l}\setminus\{\zero\}$, where $n\ge2$, $1\le d\le n$, $l=0,1$ and $dn-d^2-l\ge n-1$, then
\begin{equation}\label{2terms}
    \rho(\sjmx)=x^{(n-1)}(t_0)f_1(\sj_{t_0}^{n-3}x)+f_2(\sj_{t_0}^{n-2}x),
\end{equation}
where $f_1\in\scp_{n-2,d-1,(d-1)(n-2)-(d-1)^2-l}\setminus\{\zero\}$ and $f_2$ is a polynomial in $\sj_{t_0}^{n-2}x$. The same conclusion holds if $n\ge3$, $3\le d\le n$, $l=2$ and $dn-d^2-2\ge n-1$.
\end{lemma}
\begin{proof} If $d=1$, then the condition $dn-d^2-l\ge n-1$ forces $l=0$. Recall that $\scp_{n,1,n-1}$ is spanned by $x^{(n-1)}(t_0)$, so the conclusion of the lemma holds. If $d=2$, the conclusion of the lemma follows from \cite[Theorem 1.1]{wz}. Let $x^{(k_0)}(t_0)$ be the highest derivative on which $\rho$ depends, and let $m_0$ be the degree of $\rho$ as a polynomial in $x^{(k_0)}(t_0)$.

If $3\le d\le n$, then the conditions $\sd(\rho)=dn-d^2-l$ and $l\le2$ imply $k_0\ge n-d$. If $3\le d\le n/2$, then $k_0\ge n/2$. For any $\rho_1\in\scp_{n,d,l}$ with $d>n/2$, there exists $\rho_2\in\scp_{n,n-d,l}$ such that $\rho_1\sim_n\rho_2$ (Proposition~\ref{pnd}). Note that the highest derivative on which $\rho_2$ depends is the same as that for $\rho_1$ (see (\ref{p=1})), and the condition $dn-d^2-l\ge n-1$ is invariant under replacing $d$ with $n-d$. Hence $k_0\ge n/2$ also holds when $d>n/2$. The inequality (\ref{highest}) then yields $m_0=1$.
Therefore we may write
\begin{equation}\label{rhot}
    \rho(\sjmx)=x^{(k_0)}(t_0)f_1(\sj_{t_0}^{k_0-1}x)+f_2(\sj_{t_0}^{k_0-1}x),
\end{equation}
where $f_1$ and $f_2$ are polynomials in $\sj_{t_0}^{k_0-1}x$ and $f_1\neq\zero$. Setting $m=1$ and $i_1=k_0$ in Lemma~\ref{estimate}, we obtain
\begin{equation}\label{tf1}
    f_1\in\scp_{n-2,d-1,(dn-d^2-l)-k_0}\setminus\{\zero\}.
\end{equation}
Consequently, $f_1$ is independent of $x^{(n-1)}(t_0)$ and $x^{(n-2)}(t_0)$. It follows from Theorem~\ref{arb} that
\begin{eqnarray}
    & \sd(f_1)=(dn-d^2-l)-k_0\le(d-1)(n-2)-(d-1)^2. \hspace{1.5mm} \text{Thus} & \notag \\
    & k_0\ge n-1-l.  & \label{k0}
\end{eqnarray}

  It suffices to prove that $k_0=n-1$; in this case, the conclusion of the lemma follows from (\ref{rhot}) and (\ref{tf1}). Assume that $k_0=n-1$ for $d=d_0\ge2$ and $l=0,1$. We shall show that the same conclusion remains valid for $d=d_0+1\ge3$ and $l=0,1,2$. If $l=0$, then $k_0=n-1$ follows immediately from (\ref{k0}).

Assume for contradiction that $l=1$ and $k_0=n-2$. The conditions $dn-d^2-l\ge n-1$, $d\ge3$ and $l=1$ imply $n\ge5$, so $2k_0>n-1$. In view of (\ref{tf1}), we have $f_1\in\scp_{n-2,d-1,(d-1)(n-2)-(d-1)^2}\setminus\{\zero\}$. By the induction hypothesis, $f_1$ depends on $x^{(n-3)}(t_0)$. Hence at least one of the monomials in $\rho$ is divisible by $x^{(n-2)}(t_0)x^{(n-3)}(t_0)$, contradicting Lemma~\ref{2factors}. Therefore $k_0=n-1$ when $l=1$.

If $l=2$ and $k_0=n-2$, then $f_1\in\scp_{n-2,d-1,(d-1)(n-2)-(d-1)^2-1}$, and arguments analogous to those for $l=1$ yield a contradiction. If $l=2$ and $k_0=n-3$, then $f_1\in\scp_{n-2,d-1,(d-1)(n-2)-(d-1)^2}$, so $f_1$ depends on $x^{(n-3)}(t_0)$. Hence at least one of the monomials in $\rho$ is divisible by $x^{(n-3)}(t_0)^2$, contradicting $m_0=1$. Therefore $k_0=n-1$ in all cases.
\end{proof}

\begin{cor} \label{2}
 If $\dim\scp_{n,d,dn-d^2}\neq0$ and $\rho_0\in\scp_{n,d,dn-d^2}\setminus\{\zero\}$, where $n\ge3$ and $2\le d\le(n+1)/2$, then a constant multiple of $\rho_0$ contains the monomials
\begin{align*}
    &x^{(n-1)}(t_0)\,x^{(n-3)}(t_0)\cdots x^{(n-2d+3)}(t_0)\,x^{(n-2d+1)}(t_0)\\
    &-\frac{n-2d+3}{n-2d+2}\,x^{(n-1)}(t_0)\,x^{(n-3)}(t_0)\cdots x^{(n-2d+5)}(t_0)\,x^{(n-2d+2)}(t_0)^2
\end{align*}
(if $d=2$, the second monomial above is simply $-\frac{n-1}{n-2}x^{(n-2)}(t_0)^2$).
\end{cor}
\begin{proof}
The case $d=2$ follows from \cite[Theorem 1.1]{wz}. Assume that the conclusion of the corollary holds for $d=d_0\ge2$. We shall show that it remains valid for $d=d_0+1\ge3$. The condition $3\le d\le(n+1)/2$ implies
\[
    dn-d^2=\sum_{j=1}^{d}(n-2j+1)\ge n-1.
\]
It follows from Lemma~\ref{n-1} that $\rho_0$ has the form (\ref{2terms}). Finally we apply the induction hypothesis to the polynomial $f_1\in\scp_{n-2,d-1,(d-1)(n-2)-(d-1)^2}\setminus\{\zero\}$ in (\ref{2terms}).
\end{proof}

Let $I_d=\{(r_1,\dots,r_d)\in\bz^d: r_1,\dots,r_d\ge0\}$. We define a total order on $I_d$ as follows. For distinct elements $(r_1,\dots,r_d)$ and $(j_1,\dots,j_d)$ of $I_d$, we write
\begin{equation}\label{to}
    (r_1,\dots,r_d)\ge(j_1,\dots,j_d)
\end{equation}
if $r_1>j_1$, or if there exists $m\in\{2,\dots,d\}$ such that $r_m>j_m$ and $r_i=j_i$ for $1\le i\le m-1$.

For $\rho\in\pnd$ with $1 \le d \le n$,
we write every monomial in $\rho$ uniquely in the form $c_{(r_1,\dots,r_d)}x^{(r_1)}(t_0)\cdots x^{(r_d)}(t_0)$, where $r_1\ge\cdots\ge r_d\ge0$ and $c_{(r_1,\dots,r_d)}\in\bc^{\ast}$. The total order on $I_d$ induces a total order on the monomials in $\rho$: we say that
\[
    c_{(r_1,\dots,r_d)}x^{(r_1)}(t_0)\cdots x^{(r_d)}(t_0)\ge c_{(j_1,\dots,j_d)}x^{(j_1)}(t_0)\cdots x^{(j_d)}(t_0)
\]
if (\ref{to}) holds. Recall the linear operator $\mathfrak{M}_{(k_1,\dots,k_d)}$ from (\ref{mul}),
which can be defined on the entire space $\pn$.

\begin{lemma} \label{construct}
Let $\rho\in\pndl \setminus \{\zero\}$, where $1 \le d \le n$ and $0 \le l \le dn-d^2$, and let $c_{(r_1^0,\dots,r_d^0)}x^{(r_1^0)}(t_0)\cdots x^{(r_d^0)}(t_0)$,
where $c_{(r_1^0,\dots,r_d^0)} \in \bc^{\ast}$ and $r_1^0\ge\cdots\ge r_d^0\ge0$, be the largest monomial in $\rho$ with respect to the total order defined above. Consider the set
\[
    K=K_{d,r_d^0}=\{(k_1,\dots,k_d)\in\bz^d: 0\le k_1\le\cdots\le k_d\le r_d^0\}.
\]
Then the polynomials $\mathfrak{M}_{(k_d,k_{d-1}, \dots,k_1)}(\rho)\in\scp_{n,d,l-k_1-\cdots-k_d}$, $(k_1,\dots,k_d)\in K$, are linearly independent.
\end{lemma}
\begin{proof}
For each $(k_1,\dots,k_d)\in K$, the largest monomial in $\mathfrak{M}_{(k_d,k_{d-1}, \dots,k_1)}(\rho)$ has the form
\begin{equation}\label{large}
    c_{(r_1^0-k_1,\dots,r_d^0-k_d)}\,y^{(r_1^0-k_1)}(t_0)\cdots y^{(r_d^0-k_d)}(t_0),
\end{equation}
where $c_{(r_1^0-k_1,\dots,r_d^0-k_d)}\in\bc^{\ast}$ and $y \in \cisoc$. So $\mathfrak{M}_{(k_d,k_{d-1}, \dots,k_1)}(\rho)\neq\zero$.

Consider any finite linear combination
\begin{equation}\label{lc}
    \sum_{i=1}^{p}\lambda_i\,\mathfrak{M}_{(k_{d,i}, k_{d-1,i}, \dots,k_{1,i})}(\rho),
\end{equation}
where $\lambda_i\in\bc^{\ast}$ and $(k_{1,i},\dots,k_{d,i})\in K$ for $i=1,\dots,p$. Let $(k_{1,i_0},\dots,k_{d,i_0})$ be the smallest element of the subset
\[
    \{(k_{1,i},\dots,k_{d,i}): i=1,\dots,p\}\subset K.
\]
Then the largest monomial in (\ref{lc}) is of the form (\ref{large}) with $(k_1,\dots,k_d)=(k_{1,i_0},\dots,k_{d,i_0})$. Hence the linear combination in (\ref{lc}) is not $\zero \in \pn$.
\end{proof}

\noindent {\it Proof of Theorem \ref{bounds}. } (a) We already know that statement (a) holds when $1\le n\le5$, or when $n\ge6$ and $d=0,1,2,n-2,n-1,n$ (see Section~\ref{intro}). Suppose that (a) is valid for $2\le d=d_0\le n-4$.
We shall show that (a) remains true for $d=d_0+1$.

Let $I$ be the closed interval
\[
    \left[\frac{n-\sqrt{n^2-4n+4-4l}}{2},\ \frac{n+\sqrt{n^2-4n+4-4l}}{2}\right]\subset\br.
\]
It is straightforward to verify that $dn-d^2-l\ge n-1$ if and only if $d\in I$. The condition $[3,n-3]\subset I$ is equivalent to $l\le2n-8$. Since $n\ge6$, $l\le2$ and $3\le d\le n-3$, we have $dn-d^2-l\ge n-1$.
Statement (a) follows from Lemma~\ref{n-1} upon applying the induction hypothesis to the polynomial $f_1$ in (\ref{2terms}).

 (b) In Lemma~\ref{construct}, setting $l=dn-d^2$ and
 choosing $\rho$ to be the polynomial $\rho_0$ in Corollary~\ref{2}, we immediately obtain the first part of statement (b).
 If $l\le n-2d+1$, it follows from Proposition~\ref{monoj}
 that the monomials $x^{(k_1)}\cdots x^{(k_d)}$, where $(k_1,\dots,k_d)\in S_{n,d,l}$, form a basis for $\pndl$.
\qed



\begin{thebibliography}{HMY}
\bibitem[CL]{cl} A.L. Carey, E. Langmann,
Loop groups and quantum fields, Geometric analysis and applications to quantum field theory (Adelaide, 1998/1999), 45-94,
Progr. Math. 205,
Birkh\"{a}user Boston Inc., Boston, 2002.

\bibitem[F]{f}
E. Frenkel,
Langlands Correspondence for Loop Groups,
Cambridge Studies in Advanced Mathematics 103,
Cambridge University Press, Cambridge, 2007.

\bibitem[G]{g} M.A. Guest,
Harmonic maps, loop groups, and integrable systems,
London Math. Soc. Stud. Texts 38,
Cambridge University Press, Cambridge, 1997.

\bibitem[HSW]{hsw} N.J. Hitchin, G.B. Segal, R.S. Ward,
Integrable systems,
twistors, loop groups, and Riemann surfaces,
Lectures from the instructional conference held at the University of Oxford, Oxford, 1997.

\bibitem[HMY]{hu} H.N. Huang, S.A.M. Marcantognini,
N.J. Young, Chain Rules for Higher Derivatives, Math. Intelligencer \textbf{28} (2006), 61-69.

\bibitem[LZ]{lz} L. Lempert, N. Zhang, Dolbeault cohomology of a
loop space, Acta Math. \textbf{193} (2004), 241-268.

\bibitem[PS]{ps} A. Pressley, G. Segal, Loop groups, Oxford
University Press, New York, 1986.

\bibitem[WZ]{wz} X. Wang, N. Zhang, Holomorphic sections of line bundles on
the jet spaces of the Riemann sphere, Anal. Math. \textbf{48} (2022), 589-600.

\bibitem[Z03]{z03} N. Zhang, Holomorphic line bundles on the loop
space of the Riemann sphere, J. Differential Geom. \textbf{65} (2003), 1-17.

\bibitem[Z17]{z17} N. Zhang, Holomorphic automorphisms of the loop space of
$\mathbb{P}^n$, Comm. Anal. Geom. \textbf{25} (2017), 709-718.
\end{thebibliography}
\end{document}